\documentclass{amsart} 
\usepackage{hyperref}

\usepackage{amssymb,latexsym}
\usepackage{amsmath}
\usepackage[english]{babel}

\newcommand{\M}{\mathbb{M}}

\newcommand{\R}{\mathbb{R}}
\newcommand{\Ss}{\mathbb{S}}
\newcommand{\HH}{\mathbb{H}}

\newcommand{\sd}{s_{\delta}}
\newcommand{\cd}{c_{\delta}}

\newcommand{\trace}{\mathrm{tr\,}}
\newcommand{\Ric}{\mathrm{Ric}}
\newcommand{\iid}{\mathrm{Id}\,}

\newcommand{\lgra}{\longrightarrow}
\newcommand{\lto}{\ensuremath{\longrightarrow}}
\newcommand{\function}[5]
{\begin{eqnarray*}\begin{array}{r@{}ccl}
 #1\;\colon\;  & #2 &\lto & #3 \\[.05cm]  
  & #4 &\longmapsto  & #5 
\end{array}\end{eqnarray*}
}
\newcommand{\grad}{{\rm grad}}
\newcommand{\diff}{\mathrm{d}}

\newtheorem{thm}{Theorem}[section]

\newtheorem{lemma}[thm]{Lemma}

\newtheorem{exabout:ample}[thm]{Exemple}

\newcommand{\be}{\begin{enumerate}}  \newcommand{\ee}{\end{enumerate}}
\newcommand{\beqt}{\begin{equation}}  \newcommand{\eeqt}{\end{equation}}
\newcommand{\beq}{\begin{eqnarray}}  \newcommand{\eeq}{\end{eqnarray}}
\newcommand{\beQ}{\begin{eqnarray*}} \newcommand{\eeQ}{\end{eqnarray*}}
\begin{document}
\title{On compact embbeded Weingarten hypersurfaces in warped products}
\author[J. ROTH]{Julien Roth}
\address[J. ROTH]{Universit\'e Gustave Eiffel, CNRS, LAMA UMR 8050, F-77447 Marne-la-Vallée, France}
\email{julien.roth@univ-eiffel.fr}
\author[A. UPADHYAY]{Abhitosh Upadhyay}
 \address[A. UPADHYAY]{School of Mathematics and Computer Science, Indian Institute of Technology, Goa 403401, India}
\email{abhitosh@iitgoa.ac.in}

\maketitle

\begin{abstract}
We show that compact embedded starshaped $r$-convex hypersurfaces of certain warped products satisfying $H_r=aH+b$ with $a\geqslant 0$, $b>0$,  where $H$ and $H_r$ are respectively the mean curvature and $r$-th mean curvature is a slice. In the case of space forms,  we show that without the assumption of starshapedness, such Weingarten hypersurfaces are geodesic spheres. Finally, we prove that, in the case of space forms, if $Hr-aH-b$ is close to $0$ then the hypersurface is close to geodesic sphere  for the Hausdorff distance. We also prove an anisotropic version of this stability result in the Euclidean space.  \end{abstract} 

\section{Introduction}
The well-known Alexandrov theorem \cite{Al} ensures that a closed embedded hypersurface of the Euclidean space $\R^{n+1}$ with constant mean curvature must be a round sphere. The hypothesis, that the hypersurface is to be embedded, is crucial as proved by the counter examples of Wente \cite{We}, Kapouleas \cite{Ka} or Hsiang-Teng-Yu \cite{HTY} for instance. Further, this result has been extended to scalar curvature and then higher order mean curvatures by Ros \cite{Ro1,Ro2} as well as for any concave function of the principal curvatures by Korevaar \cite{Ko}. For higher order mean curvatures, the necessity of the embedding is still an open question.\\
\indent
Note that Montiel and Ros \cite{MR} proved that the Alexandrov theorem for the mean curvature as well as for higher order mean curvatures is also true for hypersurfaces of hyperbolic spaces and half-spheres.\\
\indent
On the other hand, very recently, de Lima \cite{dL} proved a comparable result for the so called linear Weingarten hypersurfaces satisfying $H_r=aH+b$ for two real constants $a\geqslant 0$ and $b> 0$, where $H$ and $H_r$ are respectively the mean curvature and the $r$-th mean curvature of the hypersurfaces. The hypersurfaces are supposed to be embedded in this result too and $H_r$ is a positive function.\\
\indent The aim of the present note is to show that  Lima's result also holds for a large class of warped products which contains in particular the hyperbolic spaces and the half-spheres.\\
\indent Let $n\geqslant 2$ be an integer and $(M^n,g_M)$ be a compact Riemannian manifold of dimension $n$ satisfying
$$\Ric_M\geqslant (n-1)k g,$$
for some constant $k$. Moreover, let $t_0>0$ and $h:[0,t_0)\longrightarrow \R$ be a positive function satisfying the following four conditions
\begin{equation}\label{H1}\tag{H1}
h'(0)=0\ \text{and}\ h''(0)>0,\end{equation}
\begin{equation}\label{H2}\tag{H2}
h'(t)>0\ \text{for all}\ t\in(0,t_0),\end{equation}
\begin{equation}\label{H3}\tag{H3}
\text{the function}\ r\longmapsto 2\dfrac{h''(t)}{h(t)}-(n-1)\dfrac{k-h'(t)^2}{h(t)^2}\ \text{is non-decreasing on}\ (0,t_0),
\end{equation}
\begin{equation}\label{H4}\tag{H4}
\dfrac{h''(t)}{h(t)}+\dfrac{k-h'(t)^2}{h(t)^2}>0\ \text{for all}\ t\in(0,t_0).\\
\end{equation}
We consider the warped product $P$ defined by $P=[0,t_0)\times M$ endowed with the metric $g_P=dt^2\oplus h(t)g_M$.
\begin{thm}\label{thm1}
Let $n\geqslant 2$ and $r\in\{2,\cdots,n\}$ are two integers and $\Sigma$ be a closed, oriented hypersurface embedded into the warped product $(P^{n+1},g_P)$. We assume that the four conditions \eqref{H1}-\eqref{H4} are satisfied. If $\Sigma$ is star-shaped and $H_r$ is a positive function satisfying
$H_r=aH+b$ for some real constants $a\geqslant 0$ and $b>0$, then $\Sigma$ is a slice $\{t_1\}\times M$.
\end{thm}
\noindent
Note that in \cite{WX}, Wu and Xia obtained a slightly different result for another type of relation between higher order mean curvatures.\\
We also obtain a comparable result for the space forms $\M^{n+1}(\delta)$ for which the star-shapeness is not required. Here, $\M^{n+1}(\delta)$ denotes the Euclidean space $\R^{n+1}$ if $\delta=0$, the half-sphere $\Ss^{n+1}_+(\delta)$ if $\delta>0$ and the hyperbolic space $\HH^{n+1}(\delta)$ if $\delta<0$. Namely, we have the following.
\begin{thm}\label{thm2}
\noindent
Let $n\geqslant 2$ and $r\in\{2,\cdots,n\}$ are two integers and let $\Sigma$ be a closed, connected and embedded hypersurface of $\M^{n+1}(\delta)$. Assume that the $r$-th mean curvature $H_r$ is a positive function satisfying $H_r=aH+b$ for some real constants $a\geqslant 0$ and $b>0$. Then, $\Sigma$ is a geodesic sphere of  $\M^{n+1}(\delta)$. 
\end{thm}
\noindent
This result extend the result of Lima \cite{dL}  for real space forms of non-zero sectional curvature. Note also that when $a=0$, we recover the Alexandrov theorem of Montiel and Ros \cite{MR} for hypersurfaces with constant $r$-th mean curvature.\\\\
We also consider the stability of this new characterization of geodesic spheres in space forms, precisely, the following natural question: {\it if a closed embedded and oriented hypersurface $\Sigma$  of $\M^{n+1}(\delta)$ is almost Weingarten in the following sense, $H_r=aH+b+\varepsilon$, where $\varepsilon$ is a smooth function on $\Sigma$ which is sufficiently small, is $\Sigma$ close to a geodesic sphere ?}\\
We answer this question by the following result.

\begin{thm}\label{thm3}
Let $n\geqslant 2$ and $r\in\{2,\cdots,n\}$ are two integers. Let $\Sigma$ be a closed embedded and oriented hypersurface of $\M^{n+1}(\delta)$ bounding a domain $\Omega$. There exist three constants $\gamma$, $C$ and $\varepsilon_1$, with $\gamma$ depending only on $n$; $C$ and $\varepsilon_1$ depending on $n$, $r$, $\delta$, $\displaystyle\min_{\Sigma}H_r$, $\displaystyle\min_{\Sigma}(H_{r;n,1})$, $\|B\|_{\infty}$, $V(\Sigma)$ and $R$ so that if $\Sigma$ is almost linear Weingarten in the following sense
$$H_r=aH+b+\varepsilon,$$
where $a\geqslant0$ and $b>0$ are real constants and $\varepsilon$ is a smooth function satisfying $\|\varepsilon\|_1\leqslant \varepsilon_1$, then
$$d_H(\Sigma,S_{\rho_0})\leqslant C\|\varepsilon\|_1^{\gamma},$$
where $S_{\rho_0}$ is a geodesic sphere of a certain radius $\rho_0$ and $d_H$ is the Hausdorff distance between compact sets into $\M^{n+1}(\delta)$.
\end{thm}
Here, $B$ is the second fundamental form, $V(\Sigma)$ the volume of $\Sigma$, $R$ the extrinsic radius of $\Sigma$ and $H_{r;n,1}$ is an extrinsic quantity defined from the second fundamental form (see \eqref{defHlij} for the precise definition).
\section{Preliminaries}
\subsection{Basics about warped products and Brendle's inequality}
The classical Heinze-Karcher inequality for compact embedded hypersurfaces of the Euclidean space says that
$$ \int_{\Sigma}\frac{1}{H}dv_g\geqslant (n+1)V(\Omega),$$
where $\Sigma$ is the embedded hypersurface which bounds the compact domain $\Omega$, $H$ is the mean curvature of $\Sigma$, supposed to be positive, and $V(\Omega)$ is the volume of the domain $\Omega$.\\
In \cite{Br}, Brendle proved an analogue of the Heintze-Karcher inequality for a large class of warped products manifolds, namely, the warped products $(P,g_p)$ of the form given in the introduction and satisfying conditions \eqref{H1}, \eqref{H2} and \eqref{H3}.\\
Let $\Sigma$ be a compact embedded and oriented hypersurface of $(P^{n+1}, g_P)$, we consider the function $f=h'$ and the vector field $X=h\dfrac{\partial}{\partial t}$. In \cite{Br}, Brendle used the fact that $X$ is a conformal vector field, $\mathcal{L}_Xg=2fg$, to obtain the following Hsiung-Minkowski formula for hypersurfaces of these warped products
\begin{equation}\label{HM}
\int_{\Sigma}H\langle X,\nu\rangle dv_g=\int_{\Sigma}fdv_g,
\end{equation}
where $g$ is the induced metric on $\Sigma$ and $\nu$ is the outward normal unit vector field. Using this, he was able to prove the following extension of the Heinze-Karcher inequality: 
\begin{equation}\label{HK}
\int_{\Sigma}\frac{f}{H}dv_g\geqslant (n+1)\int_{\Omega}fdv{\overline{g}}.
\end{equation}
Moreover, if equality holds, then $\Sigma$ is umbilical. If in addition, condition \eqref{H4} is satisfied, then $\Sigma$ is a slice $\{t_1\}\times N$.
\subsection{Higher order mean curvatures and Hsiung-Minkowski formulas}
The higher order mean curvatures are extrinsic quantities defined from the second fundamental form and generalising the notion of mean curvature. Up to a normalisation constant the mean curvature $H$ is the trace of the second fundamental form $B$:
\beqt
H=\frac{1}{n}\trace(B).
\eeqt
In other words, the mean curvature is  
\beqt
H=\frac{1}{n}S_1(\kappa_1,\dots,\kappa_n),
\eeqt
where $S_1$ is the first elementary symmetric polynomial and $\kappa_1,\dots,\kappa_n$ are the principal curvatures. Higher order mean curvatures are defined in a similar way for $r\in\{1,\dots,n\}$ by
\beqt 
H_r=\frac{1}{\binom{n}{r} }S_r(\kappa_1,\cdots,\kappa_n),
\eeqt
where $S_r$ is the $r$-th elementary symmetric polynomial, that is for any $n$-tuple $(x_1,\cdots,x_n)$
\beqt
S_r(x_1,\dots,x_n)=\sum_{1\leqslant i_1<\dots<i_r \leqslant n}x_{i_1}\cdots x_{i_r}.
\eeqt
 By convention we set $H_0=1$ and $H_{n+1}=0$. Finally, for convenience we also set $H_{-1}=-\langle X,\nu\rangle$.

We also recall some classical inequalities between the $H_r$ which are well-known. First, for any $r\in\{0,\cdots,n-2\}$, we have
\beqt\label{inegHr1}
H_rH_{r+2}\leqslant H_{r+1}^2,
\eeqt 
with equality at umbilical points, cf.~\cite[p.~104]{HLP}. Moreover, cf.~ \cite{BC}, if $H_{r+1}>0,$ then $H_s>0$ for any $s\in\{0,\cdots,r\}$ and 
\beqt\label{inegHr2}
H_{r+1}^{\frac{1}{r+1}}\leqslant H_r^{\frac{1}{r}}\leqslant \cdots\leqslant H_2^{\frac12}\leqslant H.
\eeqt

In \cite{KLP}, the authors prove a general weighted Hsiung-Minkowski type formula in warped product. Namely, they prove 
\begin{eqnarray}\label{KLP1}
\int_{\Sigma}\phi(fH_{k-1}-H_k\langle X,\nu\rangle )dv_g+\dfrac{1}{k\binom{n-1}{k}}\int_{\Sigma}\phi(div_{\Sigma}T_{k-1})(\xi)dv_g\\ \nonumber=-\dfrac{1}{k\binom{n-1}{k}}\int_{\Sigma}\langle T_{k-1}(\xi),\nabla \phi\rangle dv_g
\end{eqnarray}
where $\phi$ is a smooth function on $\Sigma$ and $\xi=X^T$ is the tangential part of the conformal vector field $X$. Moreover, for $r\geqslant 2$ we have
$$ (div_{\Sigma}T_{k-1})(\xi)=-\binom{n-3}{k-2}\displaystyle\sum_{j=1}^{n-1}H_{k-1;j}\xi^jRic(e_j,\nu),$$
with $H_{k-1;j}=\sigma_{k-1}(\lambda_1,\cdots,\lambda_{j-1},\lambda_{j+1},\cdots,\lambda_n)$ and where $\sigma_k$ is the $k$-th elemntary symmetric polynomial of order $k$.
In addition, if and if conditions \eqref{H1}-\eqref{H4} and $\Sigma$ is starshaped, then 
\begin{equation}\label{genHM}
 \int_{\Sigma}fH_{r-1}dv_g\leqslant \int_{\Sigma}H_r\langle X,\nu\rangle dv_g.
\end{equation}

In the case of space forms $\M^{n+1}(\delta)$, we have the classical Hsiung-Mnkowski formulas
\begin{equation}\label{hsiung0}
\int_M\Big(H_{r}\left\langle Z,\nu\right\rangle+\cd(\rho) H_{r}\Big)dv_g=0,
\end{equation}
where $\rho(x)=d(p,x)$ is the distance function to a base point $p$ (in the sequel, $p$ will be the center of mass of $M$),  $Z$ is the position vector defined by $Z=\sd(\rho)\overline{\nabla} \rho$, where $\overline{\nabla}$ is the connection of $\M^{n+1}(\delta)$ and the functions $\cd$ and $\sd$ are defined by
$$\cd(t)=\left\lbrace \begin{array}{ll}
\cos(\sqrt{\delta}t)&\text{if}\ \delta>0\\
1&\text{if}\ \delta=0\\
\cosh(\sqrt{|\delta|}t)&\text{if}\ \delta<0
\end{array}
\right. 
\quad\text{and}\quad
\sd(t)=\left\lbrace \begin{array}{ll}
\frac{1}{\sqrt{\delta}}\sin(\sqrt{\delta}t)&\text{if}\ \delta>0\\
t&\text{if}\ \delta=0\\
\frac{1}{\sqrt{|\delta|}}\sinh(\sqrt{|\delta|}t)&\text{if}\ \delta<0.
\end{array}
\right. 
$$

\subsection{Anisotropic mean curvatures}
Let $F:\Ss^n\lgra\R_+^*$ be a smooth function satisfying the following convexity assumption
\begin{equation}\label{convexity}
A_F=(\nabla dF+F\iid_{|T_x\Ss^n})_x>0,
\end{equation}
for all $x\in\Ss^n$, where $\nabla dF$ is the Hessian of $F$ and $>0$ means positive definite in the sense of quadratic forms. Now, we consider the following map 

\function{\phi}{\Ss^n}{\R^{n+1}}{x}{F(x)x+(\grad_{|S^n}F)_x.}
The image $\mathcal{W}_F=\phi(\Ss^n)$ is called the Wulff shape of $F$ and is a smooth convex hypersurface of $\R^{n+1}$ due to condition \eqref{convexity}. It is to note that if $F=1$, then the Wulff shape is the sphere $\Ss^n$.\\
\indent
Let $X:(M^n,g)\lgra\R^{n+1}$ be an isometric immersion of $n$-dimensional closed, connected and oriented Riemannian manifold $M$ into $\R^{n+1}$. We denote by $\nu$ a normal unit vector field globally defined on $M$, that is, we have $\nu:M\lgra\Ss^n$. We set $S_F=-A_F\circ\diff\nu$, where $A_F$ is defined in \eqref{convexity}. The operator $S_F$ is called the $F$-Weingarten operator or anisotropic shape operator. In this anisotropic setting, we can define all the corresponding extrinsic quantities. The anisotropic  higher order mean curvatures $H_r^F$ are defined by
$$H^F_r=\frac{1}{\binom{n}{r}}\sigma_r(S_F),$$
where $\sigma_r(S_F)$ is the $r$-th elementary symmetric polyniomial with $n$ variables computed for anisotropic principal cruvatures $\kappa_1^F,\cdots,\kappa_n^F$.

We denote simply by $H^F$ the anisotropic mean curvature $H_1^F$. Moreover, for convenience, we set $H^F_0=1$ and $H^F_{n+1}=0$ by convention. For the Wulff shape $\mathcal{W}_F$, $\kappa_1^F=\kappa_2^F=\cdots=\kappa_n^F$ are nonzero constants. Moreover, if $\kappa_1^F=\kappa_2^F=\cdots=\kappa_n^F$, then the hypersurface has to be the Wulff shape (up to homotheties and translations).\\
Like in the anisotropic case, we have the following inequalities between higher order mean curvatures. Namely, if $H^F_{r+1}>0$ then  
  $H^F_j>0$ for all $j\in\{1,\cdots,r\}$ and 
  \begin{equation}\label{newtonHF}\left(H^F_r\right)^{\frac{1}{r}}\leqslant \left(H^F_{r-1}\right)^{\frac{1}{r-1}}\leqslant\cdots\leqslant\left(H^F_2\right)^{\frac{1}{2}}\leqslant H^F.
  \end{equation}
  Moreover, in any of these inequalities, equality occurs at a point $p$ if and only if all the anisotropic principal curvatures at $p$ are equal. Hence, equality occurs everywhere if and only if $M$ is the Wulff shape $\mathcal{W}_F$, up to translations and homotheties.\\
  Finally, we recall the anisotropic Hsiung-Minkowski formulas. For $r\in\{0,\cdots n-1\}$, we have
  \begin{equation}\label{HMHF}
  \int_M\Big(F(\nu)H_r^F+H_{r+1}^F\langle X,\nu\rangle\Big)dv_g = 0.
  \end{equation}
We also have an anisotropic analogue of the Heintze-Karcher inequality (see \cite{HL}).
If $M$ is embedded (so bounds a domain $\Omega$) and $H^F$ is everywehere positive, then the following inequality holds
\begin{equation}\label{HKan} \int_M\frac{F(\nu)}{H^F}dv_g\geqslant (n+1)V(\Omega)
\end{equation}
with equality if and only if $M$ is the Wulff shape $\mathcal{W}_F$ (up to translations and homotheties).\\
\indent
All the above mentioned results by Alexandrov and Ros have analogues for anisotropic mean curvatures with the Wulff shape replacing the sphere (see \cite{HLMG}).

\subsection{Michael-Simon extrinsic Sobolev inequality}
We conclude this section of preliminaries by recalling the extrinsic Sobolev inequality of Michael and Simon \cite{MS}. If $(\Sigma,g)$ is a closed connected and oriented hypersurface of the Euclidean space, for any $\mathcal{C}^1$ function $f$ on $M$, the following inequality holds
\begin{equation}\label{sobolev}
\left(\int_Mf^{\frac{n}{n-1}}dv_g\right)^{\frac{n-1}{n}}\leqslant K(n)\int_M\left(|\nabla f|+|Hf|\right)dv_g,
\end{equation}
where $K(n)$ is a constant that depends only on $n$. Applying this inequality for the function $f\equiv 1$, we get
\begin{equation}\label{MS1}
V(\Sigma)^{\frac{n-1}{n}}\leqslant K_{n}\int_{\Sigma}|H|dv_{g}.
\end{equation}
Now, we consider $D\subset\R^{n+1}$ be an open domain and let $N^{n+1}=(D,h)$ be a conformally flat Riemannian manifold, i.e., $h=e^{2\varphi}\widetilde{h}$
where $\widetilde{h}$ is the Euclidean metric and $\varphi\in C^{\infty}(D).$ Let $(\Sigma^{n},g)\hookrightarrow (N^{n+1},h)$ be a closed, connected, oriented and isometrically immersed hypersurface. We deduce from \eqref{MS1} that
\begin{equation}\label{MS2}
V(\Sigma)^{\frac{n-1}{n}}\leqslant c_{n,\varphi}\int_{\Sigma}|\widetilde{H}|dv_{\widetilde{g}},
\end{equation}
where $\widetilde{H}$ is the mean curvature of the isometric immersion $(\Sigma^{n},\widetilde{g})\hookrightarrow (N^{n+1},\widetilde{h})$ with $\widetilde{g}=e^{-2\varphi}g$ and $c_{n,\varphi}$ is a constant depending on $n$ and $\varphi$. Note that here, $V(\Sigma)$ is the volume of $\Sigma$ with the metric $g$ which explain the dependence of the constant $c_{n,\varphi}$ on the conformal factor $\varphi$. Thus, we deduce immediately that 
\begin{equation}\label{MS3}
V(\Sigma)^{-\frac{1}{n}}\leqslant c_{n,\varphi}\|\widetilde{H}\|_1
\end{equation}
and so
\begin{equation}\label{MS3}
V(\Sigma)^{-\frac{n+1}{n}}\leqslant c_{n,\varphi}\|\widetilde{H}\|_{n+1}^{n+1}.
\end{equation}

\section{Proof of the results}
\subsection{Proof of Theorems \ref{thm1} and \ref{thm2}}
We will prove first Theorem \ref{thm1}. The proof of Theorem \ref{thm2} differs very slightly so that we will just mention after the minor differences.\\
First, by the generalized Hsiung-Minkowski formula \eqref{genHM}, we have
$$ \int_{\Sigma}\Big(fH_{r-1}-H_r\langle X,\nu\rangle\Big) dv_g\leqslant 0.$$
Using the assumption that $M$ is a Weingarten hypersurface, that is, $H_r=aH+b$, we get
\begin{equation}\label{ineg1}
\int_{\Sigma}\Big(fH_{r-1}-aH\langle X,\nu\rangle-b\langle X,\nu\rangle\Big) dv_g\leqslant 0.
\end{equation}
From the Hsiung-Minkowski \eqref{HM} formula and the divergence theorem, we have
$$\int_{\Sigma}H\langle X,\nu\rangle dv_g=\int_{\Sigma}fdv_g$$
and 
$$\int_{\Sigma}\langle X,\nu\rangle dv_g=(n+1)\int_{\Omega}fdv_{\overline{g}},$$
respectively. Hence, \eqref{ineg1} becomes
\begin{equation}\label{ineg2}
a\int_{\Sigma}fdv_g+(n+1)b\int_{\Omega}f dv_{\overline{g}}\geqslant\int_{\Sigma}fH_{r-1}dv_g.
\end{equation}
Now, since $H_r$ is supposed to be a positive function, as a consequence of \eqref{inegHr2}, we have $H_{r-1}\geqslant H_r^{\frac{r-1}{r}}$ which after reporting into \eqref{ineg2} gives
\begin{eqnarray}\label{ineg3}
a\int_{\Sigma}fdv_g+(n+1)b\int_{\Omega}f dv_{\overline{g}}&\geqslant&\int_{\Sigma}fH_{r}^{\frac{r-1}{r}}dv_g\notag\\
&\geqslant&\int_{\Sigma}fH_{r}H_r^{\frac{-1}{r}}dv_g\notag\\
&\geqslant&\int_{\Sigma}f\frac{H_r}{H}dv_g\\
&\geqslant&a\int_{\Sigma}fdv_g+b\int_{\Sigma}\frac{f}{H}dv_g,
\end{eqnarray}
where we have used in the last two lines the facts that $H_r^{\frac1r}\leqslant H$ and $H_r=aH+b$ respectively. Now, we finish by applying the Brendle inequality
$$\int_{\Sigma}\frac{f}{H}dv_g\geqslant (n+1)\int_{\Omega}fdv{\overline{g}}$$
which gives, since $b$ is positive
$$a\int_{\Sigma}fdv_g+(n+1)b\int_{\Omega}f dv_{\overline{g}}\geqslant a\int_{\Sigma}fdv_g+(n+1)b\int_{\Omega}f dv_{\overline{g}},$$
which means that all the previous inequality are in fact equalities. In particular, equality holds in the Brendle's inequality, which implies that (since condition \eqref{H4} is assumed) Theorem \ref{thm1} is proved. \hfill$\square$
\subsection{Proof of Theorem \ref{thm3}}
The strategy of the proof consists in showing that the $L^{n+1}$-norm of $\tau$ is small (compared to $\varepsilon$) and applying the following result of \cite{RS} with $p=n+1$ where $N^{n+1}$ is either the Euclidean space, the half-sphere or the hyperbolic space.
\begin{thm}{(Roth-Scheuer \cite{RS})}\label{thmRS}
Let $D\subset\R^{n+1}$ be open and let $N^{n+1}=(D,h)$ be a conformally flat Riemannian manifold, i.e., $h=e^{2\varphi}\widetilde{h}$
where $\widetilde{h}$ is the Euclidean metric and $\varphi\in C^{\infty}(D).$ Let $\Sigma^{n}\hookrightarrow N^{n+1}$ be a closed, connected, oriented and isometrically immersed hypersurface. Let $p>n\geq 2.$
Then there exist constants $c$ and $\varepsilon_{0},$ depending on $n,$ $p,$ $V(\Sigma),$ $\|B\|_{p}$ and $\|\varphi\|_{\infty},$ as well as a constant $\alpha=\alpha(n,p),$ such that whenever there holds
$$\|\tau\|_{p}\leqslant \|\widetilde{H}\|_p \varepsilon_{0},$$
there also holds
$$d_{H}(\Sigma,S_{\rho})\leqslant \frac{c\rho}{\|\widetilde{H}\|_p^{\alpha}}\|\tau\|_{p}^{\alpha},$$
 where $S_{\rho}$ is the image of a Euclidean sphere considered as a hypersurface in $N^{n+1}$ and the Hausdorff distance is also measured with respect to the metric $h.$
\end{thm}
\noindent
First, we have
\beQ
\|\tau\|_{n+1}^{2(n+1)}&=&\left( \frac{1}{V(\Sigma)}\int_M\|\tau\|^{2(n+1)}dv_g\right)^2\\
&\leqslant&\frac{1}{V(\Sigma)^2}\left(\int_M\|\tau\|^{2n}dv_g\right)\left(\int_M\|\tau\|^2dv_g\right)
\eeQ
by the Cauchy-Schwarz inequality. From this, we deduce immediately that
\begin{equation}\label{majtau}
\|\tau\|_{n+1}^{2(n+1)}\leqslant\frac{1}{V(\Sigma)}\|B\|_{\infty}^{2n}\left(\int_M\|\tau\|^2dv_g\right).
\end{equation}
Now, we estimate $\displaystyle\int_M\|\tau\|^2dv_g$. First, we have this lemma.
\begin{lemma}\label{lemtau2}
There exists a constant positive constant $K_1=K_1(n,r,\min(H_{r;n,1}),\|B\|_{\infty})$ so that
$$\|\tau\|^2\leqslant K_1\big(H_{r-1}-H_{r}^{\frac{r-1}{r}}\big).$$
\end{lemma}
\noindent{\it Proof:}
First, as mentioned in the preliminaries section, for any $k\in\{1,\cdots,n-1\}$, we have 
$$H_k^2-H_{k+1}H_{k-1}\geqslant 0.$$
We have a more precise estimate of the positivity of this term. Namely,
\begin{equation}\label{ineqJS}
H_k^2-H_{k+1}H_{k-1}\geqslant c_n\|\tau\|^2H^2_{k+1;n,1}
\end{equation}
where $c_n$ is a constant depending only on $n$ and 
\begin{equation}\label{defHlij}H_{l;i,j}=\dfrac{\partial H_{l}}{\partial \kappa_i\partial \kappa_j}=\dfrac{1}{\binom{n}{l}}\displaystyle\sum_{\begin{array}{cc}1\leqslant i_1<\cdots<i_l\leqslant n\\ i_1,\cdots,i_l\neq i,j\end{array}}\kappa_{i_1}\cdot\cdots\cdot\kappa_{i_l}.
\end{equation}
On can find the proof in \cite{Sc} for instance. Hence, for $k=r-1$, we have
\begin{equation}\label{ineqJS2}
H_{r-1}^2-H_{r}H_{r-2}\geqslant c_n\|\tau\|^2H^2_{r;n,1},
\end{equation}
which gives, with the fact that $H_{r-2}\geqslant H_r^{\frac{r-2}{r}}$,
\begin{equation}\label{ineqJS3}
H_{r-1}^2-H_{r}^{\frac{2(r-1)}{r}}\geqslant c_n\|\tau\|^2H^2_{r;n,1}.
\end{equation}
Finally, we get
\begin{equation}\label{ineqJS3}
H_{r-1}-H_{r}^{\frac{r-1}{r}}\geqslant\dfrac{ c_n\|\tau\|^2H^2_{r;n,1}}{H_{r-1}+H_{r}^{\frac{r-1}{r}}}\geqslant \dfrac{ c_n\|\tau\|^2H^2_{r;n,1}}{2H_{r-1}}.
\end{equation}
Thus, bounding $H^2_{r;n,1}$ from below by its minimum and $H_{r-1}$ form above with $\|B\|_{\infty}$, 
we get
\begin{equation}\label{ineqJS4}
H_{r-1}-H_{r}^{\frac{r-1}{r}}\geqslant\dfrac{ c_n\|\tau\|^2\min(H_{r;n,1})^2}{2\|B\|_{\infty}^{r-1}}
\end{equation}
and finally 
\begin{equation}\label{majtauHr}
\|\tau\|^2\leqslant K_1(H_{r-1}-H_{r}^{\frac{r-1}{r}})
\end{equation}
by setting $K_1=\dfrac{2\|B\|_{\infty}^{r-1}}{ c_n\min(H_{r;n,1})^2}$.\hfill$\square$\\\\
It is to note that this lemma holds independently of the fact that $M$ is almost Weingarten. Now, using this condition of being almost Weingarten, we will bound from above $H_{r-1}-H_{r}^{\frac{r-1}{r}}$. For this, we begin from the $r$-th Hsiung-Minkowski formula $$\int_{\Sigma}\left(H_r\langle X,\nu\rangle +\cd(\rho)H_{r-1}\right)dv_g=0,$$
which becomes
\begin{equation}\label{relthm41}
a\int_{\Sigma}H\langle X,\nu\rangle dv_g+b\int_{\Sigma}\langle X,\nu\rangle dv_g+\int_{\Sigma}\varepsilon\langle X,\nu\rangle dv_g+\int_{\Sigma}\cd(\rho)H_{r-1}dv_g=0,
\end{equation}
after using the assumption that $\Sigma$ is almost Weingarten. Now, using the first Hsiung-Minkowski formula
$$\int_{\Sigma}\left(H\langle X,\nu\rangle +\cd(\rho)\right)dv_g=0,$$
and the divergence theorem, \eqref{relthm41} becomes
\begin{equation}\label{relthm42}
-a\int_{\Sigma}\cd(\rho) dv_g+(n+1)b\int_{\Omega}\sd(\rho)dv_{\overline{g}}+\int_{\Sigma}\varepsilon\langle X,\nu\rangle dv_g+\int_{\Sigma}\cd(\rho)H_{r-1}dv_g=0.
\end{equation}
On the other hand, using $H\geqslant H_{r-1}^{\frac1{r-1}}\geqslant H_r^{\frac1r}$ and the Brendle inequality
\begin{eqnarray}\label{minintHr}
\int_{\Sigma}\cd(\rho)H_{r}^{\frac{r-1}{r}}dv_g&=&\int_{\Sigma}\cd(\rho)H_rH_r^{-\frac1{}{r}}dv_g\notag\\
&=&a\int_{\Sigma}\cd(\rho)HH_r^{-\frac1r}dv_g+b\int_{\Sigma}\cd(\rho)H_r^{-\frac1r}dv_g+\int_{\Sigma}\cd(\rho)\varepsilon H_r^{-\frac1r}dv_g\notag\\
&\geqslant&a\int_{\Sigma}\cd(\rho)dv_g+b\int_{\Sigma}\dfrac{\cd(\rho)}{H}dv_g+\int_{\Sigma}\cd(\rho)\varepsilon H_r^{-\frac1r}dv_g\notag\\
&\geqslant& a\int_{\Sigma}\cd(\rho)dv_g+(n+1)b\int_{\Omega}\cd(\rho)dv_{\overline{g}}+\int_{\Sigma}\cd(\rho)\varepsilon H_r^{-\frac1r}dv_g.
\end{eqnarray}
By multiplying \eqref{majtauHr} by $\cd(\rho)$ and integrating over $\Sigma$, we get
\begin{equation}\label{majtauHr2}
\int_{\Sigma}\cd(\rho)\|\tau\|^2dv_g\leqslant K_1\int_{\Sigma}\cd(\rho)(H_{r-1}-H_{r}^{\frac{r-1}{r}})dv_g.
\end{equation}
Reporting \eqref{relthm42} and \eqref{minintHr} into \eqref{majtauHr2}, we get
\begin{equation}\label{majtauHr3}
\int_{\Sigma}\cd(\rho)\|\tau\|^2dv_g\leqslant K_1\left( -\int_{\Sigma}\varepsilon  \langle X,\nu\rangle dv_g-\int_{\Sigma}\cd(\rho)\varepsilon H_r^{-\frac1r} dv_g\right),
\end{equation}
so that we deduce
\begin{equation}\label{inegtau1}
\displaystyle\inf_{\Sigma}(\cd(\rho))\int_{\Sigma}\|\tau\|^2dv_g\leqslant K_1\left[ \dfrac{\displaystyle\sup_{\Sigma}(\cd(\rho))}{\displaystyle\inf_{\Sigma}(H_r^{\frac1r})}+\displaystyle\sup_{\Sigma}(\sd(\rho))\right]\int_{\Sigma}\varepsilon dv_g.
\end{equation}
Now, we set 
$$K_1=\left\{
\begin{array}{ll}
\dfrac{K_1}{\cd(R)}\left[ \dfrac{1}{\displaystyle\min_{\Sigma}(H_r^{\frac1r})}+\displaystyle \frac{1}{\sqrt{\delta}}\right]&\text{if}\ \delta>0\\\\
K_1\left[ \dfrac{1}{\displaystyle\min_{\Sigma}(H_r^{\frac1r})}+\displaystyle R\right]&\text{if}\ \delta=0,\\\\
K_1\left[ \dfrac{ \cd(R)}{\displaystyle\min_{\Sigma}(H_r^{\frac1r})}+\displaystyle\sd(R)\right]&\text{if}\ \delta<0,
\end{array}\right.
$$
where $R$ is the extrinsic radius of $\Sigma$. Thus, we have
\begin{equation}\label{inegtau3}
\int_{\Sigma}\|\tau\|^2dv_g\leqslant  K_2\int_{\Sigma}\varepsilon dv_g,
\end{equation}
where $K_1$ is a positive constant depending on $n$, $r$, $\delta$, $\displaystyle\min_{\Sigma}H_r$, $\displaystyle\min_{\Sigma}(H_{r;n,1})$, $\|B\|_{\infty}$ and $R$. \\
Now, combining \eqref{inegtau3} and \eqref{majtau}, we get
\begin{equation}\label{inegtau4}
\|\tau\|_{n+1}^{2(n+1)}\leqslant \dfrac{K_2\|B\|_{\infty}^{2n+1}}{V(\Sigma)}\int_{\Sigma}\varepsilon dv_g=K_3\|\varepsilon\|_1,
\end{equation}
where $K_3=K_2\|B\|_{\infty}^{2n+1}$ is a constant depending on $n$, $r$, $\delta$, $\displaystyle\min_{\Sigma}H_r$, $\displaystyle\min_{\Sigma}(H_{r;n,1})$, $\|B\|_{\infty}$, $V(\Sigma)$ and $R$.\\
In order to apply Theorem \ref{thmRS}, we need to compare the $L^{n+1}$-norms of $\tau$ and the mean curvature $\widetilde{H}$ of $\Sigma$ viewed as a hypersurface of the Euclidean space after the conformal change of metric $h=e^{2\varphi}\widetilde{h}$.\\
Now, we use \eqref{MS3} to get
\begin{equation}
1\leqslant c^2_{n,\varphi}V(\Sigma)^{\frac{2(n+1)}{n}}\|\widetilde{H}\|_{n+1}^{2(n+1)}.
\end{equation}
Hence, \eqref{inegtau4} gives
\begin{eqnarray}\label{inegtau5}
\|\tau\|_{n+1}^{2(n+1)}&\leqslant&K_3c^2_{n,\varphi}V(\Sigma)^{\frac{2n+2}{n}}\|\widetilde{H}\|_{n+1}^{2(n+1)}\|\varepsilon\|_1=K_4\|\widetilde{H}\|_{n+1}^{2(n+1)}\|\varepsilon\|_1,
\end{eqnarray}
where $K_4$ is a constant depending on $n$, $r$, $\delta$, $\displaystyle\min_{\Sigma}H_r$, $\displaystyle\min_{\Sigma}(H_{r;n,1})$, $\|B\|_{\infty}$, $V(\Sigma)$ and $R$. Note that $K_4$ depends also on $\|\varphi\|_{\infty,\Omega}$ due to \eqref{MS3}, but since $\varphi$ is the conformal change of metric between $\R^{n+1}$ and $\HH^{n+1}$ or $\Ss_+^{n+1}$, this dependence can be replaced by a dependence on $\delta$ and $R$.\\
Now, if $\|\varepsilon\|_1$ is supposed to be smaller than $\varepsilon_1=\dfrac{\varepsilon_0^{2(n+1)}}{K_4}$, where $\varepsilon_0$ is the constant of Theorem \ref{thmRS}, then we have
$$\|\tau\|_{n+1}\leqslant\|\widetilde{H}\|_{n+1}\varepsilon_0,$$
so that we can apply Theorem \ref{thmRS}. Note that $\varepsilon_1$ is a positive constant depending on $n$, $r$, $\delta$, $\displaystyle\min_{\Sigma}H_r$, $\displaystyle\min_{\Sigma}(H_{r;n,1})$, $\|B\|_{\infty}$, $V(\Sigma)$ and $R$. Thus, there exists $\rho_0>0$ so that
\begin{equation}\label{inegdH1}
d_{H}(\Sigma,S_{\rho_0})\leqslant \frac{c\rho_0}{\|\widetilde{H}\|_{n+1}^{\alpha}}\|\tau\|_{n+1}^{\alpha}.
\end{equation}
Using \eqref{inegtau4} once again, we get
\begin{equation}\label{inegdH2}
d_{H}(\Sigma,S_{\rho_0})\leqslant c\rho_0 K_4^{\frac{\alpha}{2(n+1)}}\|\varepsilon\|_1^{\frac{\alpha}{2(n+1)}}=C\|\varepsilon\|_1^{\gamma},
\end{equation}
where $C=c\rho_0 K_4^{\frac{\alpha}{2(n+1)}}$ is a positive constant depending on $n$, $r$, $\delta$, $\displaystyle\min_{\Sigma}H_r$, $\displaystyle\min_{\Sigma}(H_{r;n,1})$, $\|B\|_{\infty}$, $V(\Sigma)$ and $R$ and $\gamma$ is a positive constant depending only on $n$. This concludes the proof of Theorem \ref{thm3}. \hfill$\square$

\section{An anisotropic result}
We finish this paper by an anisotropic version of Theorem \ref{thm3}.
In \cite{RU}, we prove a new characterization of the Wulff shape which is an anisotropic version of the result of de Lima for linear Weingarten hypersurfaces. Namely, we proved the following.
\begin{thm}{(Roth-Upadhyay \cite{RU})}\label{thmAdM1}
Let $n\geqslant2$ be an integer, $F:\Ss^n\lgra\R_+^*$ a smooth function satisfying the convexity assumption \eqref{convexity} and let $M$ be a closed, connected and embedded hypersurface of $\R^{n+1}$. Assume that the higher order anisotropic mean curvature $H_r^F$, $r\in\{2,\cdots n\}$ never vanishes and satisfies $H_r^F=aH^F+b$ for some real constants $a\geqslant 0$ and $b>0$. Then, up to translations and homotheties, $M$ is the Wulff shape $\mathcal{W}_F$.
\end{thm}
Also in \cite{RU}, we proved the following stability result for $r=2$.
\begin{thm}{(Roth-Upadhyay \cite{RU})}\label{thmAdM2}
Let $n\geqslant2$ be an integer, $F:\Ss^n\lgra\R_+^*$ a smooth function satisfying the convexity assumption \eqref{convexity} and let $M$ be a closed, connected and embedded hypersurface of $\R^{n+1}$. Assume that the $r$-th order anisotropic mean curvature $H_r^F$ never vanishes and satisfies $H_2^F=aH^F+b+\varepsilon$ for some real constants $a\geqslant 0$, $b>0$ and $\varepsilon$ a smooth function. Set $\rho=\left( \frac{V(M)}{V(\mathcal{W}_F)}\right)^{\frac{1}{n}}$. Then there exist a smooth parametrisation $\psi:\mathcal{W}_{\rho F}\longrightarrow M$, a vector $c_0\in\R^{n+1}$ and an explicit constant $K$ depending on $n$, $F$, $R$, $\|H_F\|_{\infty}$, $V(M)$ and $\displaystyle\inf_{\Sigma}(H_2^F)$ so that 
$$\|\psi-{\rm Id}-c_0\|_{W^{2,2}(\mathcal{W}_{\rho F})}\leqslant K \|\varepsilon\|_2.$$
\end{thm}
By comparable arguments as those used in the proof of Theorem \ref{thm3}, we can extend this result for any $r\in\{2,\cdots,r\}$. Namely, we have the following.
\begin{thm}\label{thmAdM3}
Let $n\geqslant2$ be an integer, $F:\Ss^n\lgra\R_+^*$ a smooth function satisfying the convexity assumption \eqref{convexity} and let $M$ be a closed, connected and embedded hypersurface of $\R^{n+1}$. Assume that the $r$-th anisotropic mean curvature $H_r^F$ is positive and satisfies $H_r^F=aH^F+b+\varepsilon$ for some real constants $a\geqslant 0$, $b>0$ and $\varepsilon$ a smooth function. Set $\rho=\left( \frac{V(M)}{V(\mathcal{W}_F)}\right)^{\frac{1}{n}}$. Then there exist a smooth parametrisation $\psi:\mathcal{W}_{\rho F}\longrightarrow M$, a vector $c_0\in\R^{n+1}$ and an explicit constant $K$ depending on $n$, $r$, $F$, $R$, $\|S_F\|_{\infty}$, $V(M)$, $\displaystyle\inf_{\Sigma}(H_r^F)$ and $\displaystyle\inf_{\Sigma}(H_{r;n,1}^F)$ so that 
$$\|\psi-{\rm Id}-c_0\|_{W^{2,2}(\mathcal{W}_{\rho F})}\leqslant K \|\varepsilon\|_2.$$
\end{thm}
\noindent
{\it Proof:} By computations analogous to those of Theorem \ref{thm3}, we get first that there exists a constant $A_2$  depending on $n$, $r$, $F$, $\|S_F\|_{\infty}$ and $\displaystyle\inf_{\Sigma}(H_{r;n,1}^F)$ so that
\begin{equation}\label{majtauHrF}
H_{r-1}^F-(H_{r}^F)^{\frac{r-1}{r}}\geqslant A_2\|\tau_F\|^2.
\end{equation}
Always proceeding as in Theorem \ref{thm3}, we get from the assumption that $M$ is almost anisotropic Weingarten, and using the anisotrpoic version of both Hsiung-Minkowski formula \eqref{HMHF} and Heinzte-Karcher inequality \eqref{HKan}
\begin{equation}\label{majtauHr3F}
A_2\int_{\Sigma}F(\nu)\|\tau_F\|^2dv_g\leqslant \int_{\Sigma}\varepsilon \langle X,\nu\rangle dv_g-\int_{\Sigma}F(\nu)\varepsilon H_r^{-\frac1r} dv_g.
\end{equation}
We deduce immediately from this that 
$$\|\tau_F\|_2^2\leqslant A_3\|\varepsilon\|_,$$
where $A_3$ is a constant dependin on $n$, $r$, $F$, $R$, $\|S_F\|_{\infty}$, $\displaystyle\inf_{\Sigma}(H_r^F)$ and $\displaystyle\inf_{\Sigma}(H_{r;n,1}^F)$. Note that the extrinsic radius $R$ appear here since we need to bound from above the term  $\langle X,\nu\rangle$.\\
Finally, we conclude by applying the following result of De Rosa and Gioffr\`e.
\begin{thm}[De Rosa-Gioffr\`e \cite{dRG,dRG2}]\label{thrm2}
Let $n>2$, $p\in(1,+\infty)$ and $F:\Ss^n\lgra\R_+^*$ satisfying the convexity assumption \eqref{convexity}. There exist a constant $\delta_0=\delta_0(n,p,F)>0$  such that if $M$ is closed hypersurface into $\R^{n+1}$ satisfying 
$$Vol(M)=V(\mathcal{W}_F)\quad \text{and} \quad\int_M\|\tau_F\|^pdv_g\leqslant\delta_0$$ then there exist a smooth parametrisation $\psi:\mathcal{W}_F\longrightarrow M$, a vector $c_0\in\R^{n+1}$ and a constant $C$ depending on $n,p$ and $F$ so that 
$$\|\psi-{\rm Id}-c_0\|_{W^{2,p}(\mathcal{W}_F)}\leqslant C\|\tau_F\|_p.$$
Moreover, if $p\in(1,n]$, then the condition $\int_M\|\tau_F\|^pdv_g\leqslant\delta_0$ can be dropped.
\end{thm}
Here, it is important to mention that the volume of $M$ is supposed to be equal to $V(\mathcal{W}_F)$. If we do not assume this, the same holds replacing $\mathcal{W}_F$ by the homothetic of $\mathcal{W}_F$ of volume equal to $V(M)$, that is for $\mathcal{W}_{\rho F}$ for $\rho=\left( \frac{V(M)}{V(\mathcal{W}_F)}\right)^{\frac{1}{n}}$. Since in the statement of Theorem \ref{thmAdM3}, we do not assume that the volume is equal to $V(\mathcal{W}_F)$, this introduce a dependence of the constant $C$ also on $V(M)$. Note that this is also the case for Theorem \ref{thm3}.
\hfill$\square$
\section*{\textbf{Acknowledgements}}
Second author gratefully acknowledges the financial support from the Indian Institute of Technology Goa through Start-up Grant  (\textbf{2021/SG/AU/043}).


\begin{thebibliography}{00}
\bibitem {Al} A.D. Alexandrov, {\em A characteristic property of spheres}, Ann. Mat. Pura Appl., {\bf 58} (1962), 303--315.
\bibitem{BC} J. Barbosa and A. {\em Colares, Stability of hypersurfaces with constant r-mean curvature}, Ann. Glob. Anal. Geom. {\bf 15} (1997), no. 3, 277–297.
\bibitem{Br} S. Brendle, {\it Constant mean curvature surfaces in warped product manifolds}, Publ. Math. IHES, {\bf 117} (2013) 247-269.
\bibitem{dL} E. De Lima, \emph{A note on compact Weingarten hypersurfaces embedded in $\R^{n+1}$}, Arch. Math. (Basel), {\bf 111}, (2018), 669--672.
\bibitem{dRG} A. De Rosa \& S. Gioffr\`e, \emph{Quantitative stability for anisotropic nearly umbilical hypersurfaces}, J. Geom. Anal., {\bf 29} (2019) 2318-2346.
\bibitem{dRG2} A. De Rosa \& S. Gioffr\`e, \emph{Absence of bubbling phenomena for non convex anisotropic nearly umbilical and quasi Einstein hypersurfaces}, J. Reine Angew. Math., {\bf 780} (2021), 1-40.
\bibitem{HLP} G. Hardy, J. Littlewood \& G. Polya, {\em Inequalities}, Cambridge University Press, 1952.
\bibitem{HL} Y. He \& H. Li, \emph{Integral formula of Minkowski type and new characterization of the Wulff shape}, Acta Math. Sinica, {\bf 24(4)} (2008), 697-704.
\bibitem{HLMG} Y. He, H. Li, H. Ma \& J. Ge, \emph{Compact embedded hypersurfaces with constant higher order anisotropic mean curvature}, Indiana Univ. Math. J. {\bf 58} (2009), 853--868.
\bibitem{Ho} H. Hopf, {\em Differential Geometry in the large}, Lecture Notes in Mathematics, 1000, Springer-Verlag (1983).
\bibitem{HTY} W. Hsiang, Z. Teng \& W. Yu, {\em New examples of constant mean curvature immersions of (2k-1)-spheres into Euclidean 2k-space}, Ann. Math., {\bf 117(3)} (1983), 609--625 .
\bibitem{Hs} C. C. Hsiung, \emph{Some integral formulas for closed hypersurfaces}, Math. Scand. {\bf 2} (1954), 286-294.
\bibitem{Ka} N. Kapouleas, {\em Constant mean curvature surfaces constructed by fusing Wente tori}, Invent. Math. , {\bf 119(3)}, (1995), 443--518.
\bibitem{Ko} N. J. Korevaar, \emph{ Sphere theorems via Alexandrov constant Weingarten curvature hypersurfaces: appendix to a note of A. Ros}, J. Differential Geom. {\bf 27} (1988), 221-223.
\bibitem{KLP} K. K. Kwong, H. Lee \& J. Pyo, \emph{Weighted Hsiung-Minkowski formulas and rigidity of umbilical hypersurfaces}, Math. Res. Letters, 25(2), 597-616.
\bibitem{MS} Michael \& Simon, {\em Sobolev and mean-value inequalities on generalized submanifolds of $\R^n$}, Comm. Pure Appl. Math. {\bf 26} (1973), no. 3, 361-379.
\bibitem{MR} S. Montiel \& A. Ros, {\em Compact hypersurfaces : the Alexandrov theorem for higher order mean curvature}, Pitman Monographs Surveys Pure Appl. Math. {\bf 52} (1991), 279-296.
\bibitem{Ro1} A. Ros, \emph{Compact hypersurfaces with constant scalar curvature and a congruence theorem} J. Differential Geom. {\bf 27(2)} (1988), 215-223. 
\bibitem{Ro2} A. Ros, \emph{Compact hypersurfaces with constant higher order mean curvatures}, Rev. Mat. Iberoamericana {\bf 3(3-4)} (1987), 447-453.
\bibitem{RS} J. Roth \& J. Scheuer, \emph{Explicit rigidity of almost-umbilical hypersurfaces}, Asian J. Math. {\bf 22 (6)} (2018), 1075-1088.
\bibitem{RU} J. Roth \& A. Upadhyay, \emph{On compact anisotropic Weingarten hypersurfaces in Euclidean space},  Arch. Math (Basel) {\bf 113} (2019), no 2,  213-224.
\bibitem{Sc} J. Scheuer, \emph{Stability from rigidity via umbilicity}, preprint,  arxiv:2103.07178.
\bibitem{We} H. Wente, {\em Counterexample to a conjecture of {H}. {H}opf}, Pac. J. Math., {\bf 121(1)} (1986), 
193--243.
\bibitem{WX} J. Wu \& C. Xia, {\em On rigidity of hypersurfaces with constant curvature functions in warped product manifolds}, Ann. Glob. Anal. Geom. {\bf 46} (2014), 1-22.

\end{thebibliography}
\end{document}